\documentclass[letter, 10pt]{amsart}

\usepackage{amsthm}
\usepackage{amsmath,latexsym}
\usepackage{amssymb}
\usepackage{color}
\usepackage{setspace}
\usepackage[english, activeacute]{babel}
\usepackage[latin1]{inputenc}
\usepackage{float}
\usepackage{hyperref}
\usepackage{graphicx}
\usepackage{pdflscape}

\usepackage[left=2.7cm, right=2.7cm, top=2.5cm, bottom=2.7cm]{geometry}

\newtheorem{theo}{Theorem}[section]
\newtheorem{lemma}{Lemma}[section]
\newtheorem{corollary}{Corollary}[section]
\newtheorem{proposition}{Proposition}[section]
\newtheorem{remark}{Remark}[section]

\newcommand{\Frob}{\mathrm{F}}

\newcommand{\sgfs}{short generating functions }
\newcommand{\g}{\mathrm{g}}
\newcommand{\G}{\mathrm{G}}
\newcommand{\m}{\mathrm{m}}
\newcommand{\Ap}{\mathrm{Ap}}

\newcommand{\Z}{\mathbb{Z}}
\newcommand{\N}{\mathbb{N}}
\newcommand{\R}{\mathbb{R}}

\newcommand{\dsum}{\displaystyle\sum}
\newcommand{\dlim}{\displaystyle\lim}
\newcommand{\dprod}{\displaystyle\prod}

\title{Computing the number of numerical semigroups using generating functions}

\author[V. Blanco]{V\'ictor Blanco}
\author[P.A. Garc\'ia-S\'anchez]{Pedro A. Garc\'ia-S\'anchez}
\address{Departamento de \'Algebra, Universidad de Granada}
\email{vblanco@ugr.es, pedro@ugr.es}

\author[J. Puerto]{Justo Puerto}
\address{Departamento de Estad\'istica e Investigaci\'on Operativa, Universidad de Sevilla}
\email{puerto@us.es}
\date{\today}



\begin{document}
\begin{abstract}
This paper presents a new methodology to count the number of numerical semigroups of given genus or Frobenius number. We apply generating function tools to the
 bounded polyhedron that classifies the semigroups with given genus (or Frobenius number) and multiplicity. First, we give theoretical results
 about the polynomial-time complexity of counting the number of these semigroups. We also illustrate the methodology analyzing the cases of multiplicity $3$
and $4$ where some formulas for the number of numerical semigroups for any genus and Frobenius number are obtained.
\end{abstract}

\maketitle

\section{Introduction}
A numerical semigroup is a subset $S$ of $\N$ that is closed under
addition, $0\in S$ and generates $\Z$ as a group. This last
condition is equivalent to $\gcd(S) = 1$.

For a given numerical semigroup $S$, the set $\G(S) = \N \setminus S$, known as the set of gaps of $S$, has
finitely many elements. Furthermore, $S$ has a unique minimal
system of generators $\{n_1 < \ldots < n_p\}$. The element $n_1$ is
the least positive integer belonging to $S$ and it is denoted by
$\m(S)$, the \textit{multiplicity of} $S$ ($\m(S) = \min(S\setminus
\{0\})$) and the cardinality of $\G(S)$ is known
as the \textit{genus} of $S$, $\g(S)$. The largest integer not in $S$ is known as the Frobenius number of $S$ and it is denoted by $\Frob(S)$. The interested reader is referred
to \cite{garcia-rosales09} for further details on numerical semigroups.

Given $n \in S\setminus \{0\}$, the Ap\'ery set of $S$ with respect to $n$ is the set $\Ap(S,n) = \{s
\in S :  s - n \not\in S\}$ and it can be easily shown that if for
every $i \in \{0, \ldots, n - 1\}$ we take $w(i)$ the least
element in $S$ congruent with $i$ modulo $n$ (denoted $w(i) \equiv
i \pmod n)$, then $\Ap(S,n) = \{0 = w(0), w(1), \ldots,  w(n -
1)\}$. The set $\Ap(S,n)$ completely determines $S$, since $S =
\langle \Ap(S,n) \cup \{n\} \rangle$ (where $\langle A \rangle$
denotes the monoid generated by $A$). Moreover, the set $\Ap(S,n)$
contains, in general, more information than an arbitrary system of
generators of $S$. For instance,
Selmer in \cite{selmer77} gives the formulae,
$\g(S)=\frac{1}{\m(S)}\left(\sum_{w \in \Ap(S, \m(S))} w\right) -
\frac{m-1}{2}$ and $\Frob(S) = \max(\Ap(S,n)) - n$. Moreover, for all $s \in S$ there exist unique $t \in \N$ and $w \in \Ap(S,n)$
such that $s = tn + w$. The smallest Ap\'ery set is $\Ap(S, \m(S))$.

Rosales et al. \cite{garcia-rosales-garcia2-branco02} and previously Kunz in \cite{kunz87} gave a one-to-one correspondence between the set of
 numerical semigroups with multiplicity $m$ and the set of integer points inside a rational polyhedron in $\R^{m-1}$. This correspondence is
 based on the Ap\'ery set description of these numerical semigroups. These polyhedra are, in general, not bounded, but performing
adequeate cuts to them, like fixing the gender or the Frobenius number, we obtain a polytope whose integer points are identified
with each of the numerical semigroups with multiplity $m$ and fixed gender or Frobenius number, respectively. Then, the problem of counting numerical semigroups is equivalent to the problem
 of counting the number of integer points inside a polytope. However, it is well-known that the
problem of detecting a lattice point in polyhedra is NP-hard \cite{garey-johnson79}. In this paper,
 we apply some results on short generating functions to provide some new complexity results
 about the task of counting the number of numerical semigroups, and also some explicit descriptions for these numbers for multiplicities $3$ and $4$.


A few papers have appeared recently analyzing the number of numerical semigroups with given genus or Frobenius number \cite{bras08, bras09, rosales05}.
Here we give new results on the complexity of counting the number of numerical semigroups when the genus or the Frobenius number are fixed using \sgfs\hspace*{-0.15cm}.
A constructive approach for counting numerical semigroups and some conjectures about these numbers have been recently
presented in \cite{gap-numsem} and \cite{bras08}. In those approaches, the computation of the set of numerical semigroups with genus $g$
 requires to compute previously those with genus $g-1$, proceeding recursively. Hence, giving exponential-time algorithms in $g$
 and consuming a considerable amount of CPU memory. However, although the proposed methodology here has nice theoretical properties, from a computational viewpoint the methods given in \cite{gap-numsem} or \cite{bras08} are capable to compute all the numerical semigroups for genus up to $50$,
while we were only able to compute the number of numerical semigroups with genus up to $15$ (see Table \ref{numsem0to15}). One of the advantages of our methodology is that for given multiplicity, the polytope obtained fixing the genus (respectively the Frobenius number)
 can be seen as a parametric polytope, and a quasi-polynomial description of the number of numerical semigroups with any genus (respectively Frobenius number) can be computed.

The paper is organized as follows. Section \ref{barvinok} recalls
the main notions and results on rational generating functions. Section \ref{sec:alg} is devoted to prove
 the complexity results about counting the number of numerical semigroups, maximal embedding dimension numerical semigroups with given genus, and a
 table with the results for genus up to 15. We also provide in this section a complexity result about computing the number of numerical semigroups with given Frobenius number and multiplicity.
In Section \ref{sec:cases} we give explicit formulas for the number of numerical semigroups and maximal embedding dimension numerical semigroups with multiplicities $3$ and $4$, fixing genus and Frobenius number with multiplicities three and four.

\section{Short generating functions}
\label{barvinok}
Short generating functions were used by Barvinok
\cite{barvinok03}, initially as a tool for counting the number of
integer points inside convex polytopes. This tool is based in the
geometrical papers by Brion \cite{brion88}, Khovanskii and Puhlikov
\cite{khovanskii-puhlikov92}, and Lawrence \cite{lawrence91}.

Let $P = \{x \in \R^n: A\,x \leq b\}$ be a rational polyhedron in
$\R^n$. The main idea is to encode the
integer points inside a rational polytope as a ``long" sum of
monomials:
$$
f(P;z) = \dsum_{\alpha\in P\cap\Z^n}\,z^\alpha,
$$
where $z^\alpha = z_1^{\alpha_1}\cdots z_n^{\alpha_n}$.

Barvinok's aimed goal was representing that formal sum of monomials in
the multivariate polynomial ring $\Z[z_1, \ldots, z_n]$, as a
``short'' (polynomially indexed, for fixed $n$) sum of rational
functions.  Actually, in \cite{barvinok94}, Barvinok provides a
polynomial-time algorithm when the dimension, $n$, is fixed, to
compute those functions.
\begin{theo}[Theorem 5.4 in \cite{barvinok94}]
\label{theo:barvinok} Assume $n$, the dimension, is fixed. Given a
rational polyhedron $P \subset \R^n$ , the generating function
$f(P;z)$ can be computed in polynomial time in the form
$$
f(P;z) = \dsum_{i\in I} \varepsilon_i
\dfrac{z^{u_i}}{\dprod_{j=1}^n (1-z^{v_{ij}})}
$$
where $I$ is a polynomial-size indexing set, and where
$\varepsilon \in \{1,-1\}$ and $u_i , v_{ij} \in \Z^n$ for all $i$
and $j$.
\end{theo}
As a corollary of this result, Barvinok gave an algorithm for
counting the number of integer points in $P$. It is clear from the
original expression of $f(P;z)$ that this number is
$f(P;\mathbf{1})$, but $\mathbf{1}=(1, \ldots, 1)$ is a pole for
the rational function. Hence, the number of integer points in the
polyhedron is $\dlim_{z\rightarrow \mathbf{1}} f(S;z)$. This limit
can be computed using residue calculation tools from elementary
complex analysis.

A new algorithm for counting the integer points inside convex
polytopes using \sgfs was recently developed by Verdoolaege and Woods
\cite{woods-verdoolaege08} and implemented in the software
\texttt{barvinok} \cite{soft:barvinok}. This software also allows to
count integer points inside parametric polytopes \cite{verdoolaege07}

The above approach, apart from counting lattice points, has been
applied to solve Operations Research problems. Actually, integer
programs~\cite{deloera04, woods-yoshida05}, multiobjective integer
programs~\cite{blanco-puerto08b}, bilevel integer
programs~\cite{koeppe09} or integer programming games~\cite{koeppe08},
among many others problems are studied using Barvinok's rational
functions. Furthermore, enumerating the complete set of solutions of a
diophantine system of inequalities (in the bounded case) may be seen
as multicriteria integer problem, and then, a polynomial-delay
algorithm (for fixed dimension) have been developed in
\cite{blanco-puerto08b} for this task.
\section{Using Barvinok's algorithms to count numerical semigroups with given genus}
\label{sec:alg}
In this section we present some results concerning the task of counting the number of semigroups with fixed
multiplicity and genus or Frobenius number. These results are based on the transformation of that problem to the problem of counting the number of integer points inside certain convex polytope. This is when Barvinok's rational function theory comes into scene.

For the sake of completeness, we include the following technical result that we use in Theorem \ref{theo:numsem}.
\begin{lemma}
\label{propo1}
Let $g(x) = \dfrac{p(x)}{q(x)}$, where $p(x)$ and $q(x)$ are polynomials in the indeterminate $x$, and such that $1$ is a root of $q(x)$ with multiplicity $r$. If $\dlim_{x\rightarrow 1} g(x)$ exists, then $p^{k)}(1)=0$ for $k=1, \ldots, r$.
\end{lemma}
\begin{proof}
Note that since $1$ is a root of $q(x)$ with multiplicity $r$, we can write $q(x) = (1-x)^r h(x)$ with $h(x)$ a polynomial with $h(1)\neq 0$.

Assume that $p^{k)}(1) \neq 0$ for $k\leq r$. Clearly, $q^{k)}(1)=0$, and then, $\dlim_{x\rightarrow 1} g(x) = \dlim_{x\rightarrow 1} \dfrac{p(x)}{q(x)} = \cdots = \dlim_{x\rightarrow 1}
\dfrac{p^{k)}(x)}{q^{k)}(x)} = \dfrac{p^{k)}(1)}{0} = \infty$, a contradiction with the existence of limit of $g(x)$ at $x=1$.
\end{proof}
\begin{theo}
\label{theo:numsem}
For fixed multiplicity $m$, counting the numerical semigroups for any genus $g$ is doable in polynomial time.
\end{theo}
\begin{proof}
Let $\mathcal{S}(m)$ be the set of all numerical semigroups with
multiplicity $m \in \N\setminus \{0, 1\}$. In \cite{garcia-rosales-garcia2-branco02}, Rosales et al. proved that
there is a one-to-one correspondence between this set and the set
of non-negative integer solutions of a system of linear diophantine inequalities. This identification was previously used by Kunz in \cite{kunz87}.

Let $S$ be in $\mathcal{S}(m)$ with $\Ap(S,m)=\{0=w(0), w(1), \ldots, w(m-1)\}$. For all $i \in  \{1, \ldots, m-1\}$, let $k_i \in \N$ be such that $w(i)=k_i\,m + i$. Then $(k_1, \ldots, k_{m-1})$ is a non-negative solution of the system
\begin{eqnarray}
x_i  \geqslant&1 & \mbox{for all $i \in \{1, \ldots, m-1\}$,}\nonumber\\
x_i+x_j-x_{i+j} \geqslant&0 & \mbox{for all $1 \leqslant i \leqslant j \leqslant m-1$, $i+j \leqslant m-1$,}\label{eq:diophantine}\\
x_i+x_j-x_{i+j-m}  \geqslant&-1 &\mbox{for all $1 \leqslant i
\leqslant j \leqslant m-1$, $i+j > m$},\nonumber\\
x_i &\in \Z & \mbox{for all $i \in \{1, \ldots, m-1\}$.}\nonumber
\end{eqnarray}
Denote by $\mathcal{T}(m)$ the set of non-negative solutions of
\eqref{eq:diophantine}. Then, $\mathcal{T}(m)$ and $\mathcal{S}(m)$ are one-to-one identified (Theorem 11 in \cite{garcia-rosales-garcia2-branco02}).

Selmer in \cite{selmer77} proved that if $\Ap(S, m) = \{0, w(1)=k_1m+1, \ldots, w(n-1)=k_{m-1}+m-1\}$, $\g(S)= \frac{1}{\m(S)}\left(\sum_{w \in \Ap(S, \m(S))} w\right) - \frac{m-1}{2}$, and then then $g = \dsum_{i=1}^{m-1} k_i$.

Thus, the following system describes completely the family
$\mathcal{S}_{m,g}$, the set of numerical semigroups with multiplicity
$m$ and genus $g$. 
\begin{eqnarray}
x_i  \geqslant&1 & \mbox{for all $i \in \{1, \ldots, m-1\}$,}\nonumber\\
x_i+x_j-x_{i+j} \geqslant&0 & \mbox{for all $1 \leqslant i \leqslant j \leqslant m-1$, $i+j \leqslant m-1$,}\nonumber\\
x_i+x_j-x_{i+j-m}  \geqslant&-1 &\mbox{for all $1 \leqslant i
\leqslant j \leqslant m-1$, $i+j > m$,}\label{eq:polytope}\\
\dsum_{i=1}^{m-1} x_i &= g,&\nonumber\\
x_i \in \Z & &\mbox{for all $i \in \{1, \ldots, m-1\}$,}\nonumber
\end{eqnarray}
Furthermore, \eqref{eq:polytope} defines a polytope (bounded polyhedron) in $\R^{m-1}$, that we will call $P_{m,g}$.

Each element in $P_{m,g} \cap \Z^{m-1}$ corresponds with a numerical
semigroup with multiplicity $m$ and genus $g$. We denote $n_{m,g} = \#(P_{m,g} \cap \Z^{m-1})$.

By Theorem \ref{theo:barvinok}, the short generating function encoding the integer points inside $P_{m,g}$ is computable in polynomial time for fixed $m$ (that in this case is redundant since $m$ is clearly fixed). Let $f(P_{m,g}; z) = \dsum_{i\in I}
\varepsilon_i \dfrac{z^{u_i}}{\dprod_{j=1}^n (1-z^{v_{ij}})}$ be that generating function. Choose $c \in \Z^m$ such that $c\,v_{ij} \neq 0$ for
all $i,j$ (see \cite{barvinok03} for further details about the polynomial-time complexity of this choice). Do the changes $x_i = t^{c_i}$ for each $i=1,\ldots, m$. Then, after those changes, let $f(t):= f(P_{m,g}; t) = \dsum_{i\in I} \varepsilon_i \dfrac{t^{c\,u_i}}{\dprod_{j=1}^n (1-t^{c\,v_{ij}})}$.

Clearly, we can write $f(t)=\dfrac{P(t)}{(1-t)^r\,Q(t)}$, where $P(t)$ and $Q(t)$ are polynomial with $Q(1)\neq 0$. By using L'Hopital Rule (it can be done by Lemma \ref{propo1}) sequentially to compute the limit:
$$
\lim_{z\rightarrow 1} f(P_{m,g};z) = \lim_{t\rightarrow 1}
f(P_{m,g};t) = \lim_{t\rightarrow 1} \dfrac{P(t)}{(1-t)^r\,q(t)} = \dfrac{P^{r)}(1)}{r!q(1)} = n_{m,g}.
$$
These operations are clearly polynomially bounded, and then, the overall procedure runs in polynomial time.
\end{proof}
For a numerical semigroup with fixed genus $g$, its multiplicity is at most $g+1$ (this case is achived at the numerical semigroup $\langle g+1, \ldots, 2g+1\rangle$. Then, the number of semigroups with fixed genus $g$
(independently of the multiplicity) is given by the finite sum $n_g = \dsum_{m=2}^{g+1} n_{m,g}$. By applying Theorem \ref{theo:numsem} $g+1$ times we obtain the following result.
\begin{corollary}
\label{coro:numsem}
Let $g$ be a fixed genus, counting the number of numerical semigroups of genus $g$ is doable in polynomial time.
\end{corollary}
We have implemented the steps in the above proof to obtain the number of numerical semigroups of genus up to $15$.
Table \ref{numsem0to15} shows the number of numerical semigroups from genus $2$ to $15$ obtained using the above
methodology. We run the software \verb"barvinok"~\cite{soft:barvinok} for counting $n_{m,g}$, for each pair $(m,g)$
 with $2\leq m \leq g+1$ in a PC with an Intel Pentium 4 processor at 2.66GHz and 1 GB of RAM. It was able to compute
 the number of numerical semigroups for genus up to $15$. Although this methodology can be applied to any genus, the
software \verb"barvinok" fails to compute the short generating function representation of the polytope in \eqref{eq:polytope}
 for $g>15$. \verb"barvinok" is designed to be applicable to general polytopes and therefore does not exploit the special
structure of the polytope \eqref{eq:polytope}.
Up to date, there are only few implementations for computing generating functions of
 rational polyhedra. \texttt{barvinok} seems to be the most recent software and
incorporates a large battery of options. More effective implementations oriented to the particular polytope
 in \eqref{eq:polytope} and further research in the computation of generating functions of polyhedra would
 permit go beyond the case $g=16$. It is clear that, although
theoretically our methodology is better than the complete enumeration proposed in \cite{bras08, bras09},
 these other methods are faster than the one proposed here (the package \texttt{numericalsgps} \cite{gap-numsem} takes 0.08 seconds to compute the whole set of numerical semigroups up to genus $15$ while our algorithm takes 38.85 seconds to count the set of numerical semigroups with multiplicity $9$ and genus $12$). It seems to be due to the fact that the current
 implementation in \texttt{barvinok} does not take into account the special structure of the polytopes describing numerical semigroups.
\begin{table}[h]
{\scriptsize\begin{center}
\begin{tabular}{c|rrrrrrrrrrrrrrr|r}
{$g \backslash m$} &   2 &   3 &   4 &   5 &   6 &   7 &   8 &   9 &  10 &  11 &  12 &  13 &  14 & 15 & 16 &$n_g$ \\\hline
  2 & 1 & 1 &   &   &   &   &   &   &   &   &   &   & &&  & 2 \\
  3 & 1 & 2 & 1 &   &   &   &   &   &   &   &   &   &   &&& 4 \\
  4 & 1 & 2 & 3 & 1 &   &   &   &   &   &   &   &   &  && & 7 \\
  5 & 1 & 2 & 4 & 4 & 1 &   &   &   &   &   &   &   &   & && 12 \\
  6 & 1 & 3 & 6 & 7 & 5 &  1 &   &   &   &   &   &   &   &  &&33 \\
  7 & 1 & 3 & 7 &  10 &  11 & 6 & 1 &   &   &   &   &   &   & && 39 \\
  8 & 1 & 3 & 9 &  13 &  17 &  16 & 7 & 1 &   &   &   &   &   &  &&67 \\
  9 & 1 & 4 &  11 &  16 &  27 &  28 &  22 & 8 & 1 &   &   &   &   &&& 118 \\
 10 & 1 & 4 &  13 &  22 &  37 &  44 &  44 &  29 & 9 & 1 &   &   &   &&& 204 \\
 11 & 1 & 4 &  15 &  24 &  49 &  64 &  72 &  66 &  37 &  10 & 1 &   &   &&& 343 \\
 12 & 1 & 5 &  18 &  32 &  66 &  85 & 116 & 116 &  95 &  46 &  11 & 1 &   &&& 592 \\
 13 & 1 & 5 &  20 &  35 &  85 & 112 & 172 & 188 & 182 & 132 &  56 &  12 & 1&& &1001\\
 14 &     1     & 5     & 23    & 43    & 106   & 148   & 239   & 288   & 304   & 277   & 178   & 67    & 13    & 1 & &  1693 \\
 15 &     1     & 6     & 26    & 51    & 133   & 191   & 325   & 409   & 492   & 486   & 409   & 234   & 79    & 14    & 1 & 2857\\\hline
\end{tabular}
\end{center}}
\caption{Number of numerical semigroups with given genus $g$ and multiplicity $m$.\label{numsem0to15}}
\end{table}
Next we apply a similar methodology but fixing the Frobenius number instead of the genus, that is, our goal here is to count the number of numerical semigroups with fixed multiplicity and Frobenius number.

Let $S$ be a numerical semigroup with multiplicity $m$ and Frobenius number $F$. The key here is that the Frobenius number of $S$ can be computed by using the Ap\'ery set as $\Frob(S) = \max \Ap(S, m) - m$ (see \cite{selmer77}). Then, adding to the system of diophantine inequalities \eqref{eq:diophantine} those related to the Frobenius number:
\begin{equation*}
m\,x_i + i \leq F + m, 
\end{equation*}
and the one that fixes the element in the \'Apery set that reaches the
maximum: $m\,x_{k^*} + k^* = F+m$ (where $k^*= F \pmod m$), we have a
system of inequalities describing a polytope that characterizes
completely the set of numerical semigroups of multiplicity $m$ and
Frobenius number $F$. 
\begin{eqnarray}
x_i  \geqslant&1 & \mbox{for all $i \in \{1, \ldots, m-1\}$,}\nonumber\\
x_i+x_j-x_{i+j} \geqslant&0 & \mbox{for all $1 \leqslant i \leqslant j \leqslant m-1$, $i+j \leqslant m-1$,}\nonumber\\
x_i+x_j-x_{i+j-m}  \geqslant&-1 &\mbox{for all $1 \leqslant i
\leqslant j \leqslant m-1$, $i+j > m$,}\label{eq:polyFrob}\\
m\,x_i + i &\leq F + m &\mbox{for all $i \in \{1, \ldots, m-1\}$},\nonumber\\
m\,x_{k^*} + k^* &= F+m,&\nonumber\\
x_i \in \Z & &\mbox{for all $i \in \{1, \ldots, m-1\}$},\nonumber
\end{eqnarray}
where $k^*= F \pmod m$.

The following result follows from the above reasoning.
\begin{theo}
For fixed multiplicity $m$, counting the numerical semigroups for any Frobenius number $F$ is doable in polynomial time.
\end{theo}
\begin{remark}[Maximal Embedding dimension semigroups]
A numerical semigroup is a maximal embedding dimension numerical semigroup (MED-semigroup for short) if its multiplicity equals its embedding dimension. Further details about these semigroups can be found in \cite{barucci, rosales-med04}. In \cite{garcia-rosales-garcia2-branco02} the authors give a similar characterization that the one used in the proof of Theorem \ref{theo:numsem} but for MED-semigroups. Let $\mbox{MEDS}(m)$ be the set of all MED-semigroups with multiplicity $m \in \N\setminus \{0\}$. Then, there is a one-to-one correspondence between $\mbox{MEDS}(m)$ and the set
of solutions of the following system of linear diophantine inequalities:
\begin{eqnarray}
x_i  \geqslant&1 & \mbox{for all $i \in \{1, \ldots, m-1\}$,}\nonumber\\
x_i+x_j-x_{i+j} \geqslant&1 & \mbox{for all $1 \leqslant i \leqslant j \leqslant m-1$, $i+j \leqslant m-1$,}\label{eq:meds}\\
x_i+x_j-x_{i+j-m}  \geqslant&0 &\mbox{for all $1 \leqslant i
\leqslant j \leqslant m-1$, $i+j > m$},\nonumber\\
x_i &\in \Z & \mbox{for all $i \in \{1, \ldots, m-1\}$}.\nonumber
\end{eqnarray}
With this characterization and following the same scheme that the one used in the proof of Theorem \ref{theo:numsem} we can state that for fixed multiplicity $m$, counting the MED-semigroups for any genus $g$ is doable in polynomial time. Also, the analogous of Corollary \ref{coro:numsem} is true: Let $g$ be a fixed genus, counting the number of MED-semigroups of genus $g$ is doable in polynomial time.

Table \ref{med0to15} summarizes the results obtained up to genus 15.
\end{remark}
\begin{table}[h]
{\scriptsize\begin{center}
\begin{tabular}{c|rrrrrrrrrrrrrrr|r}
{$g \backslash m$} &	2	&	3	&	4	&	5	&	6	&	7	&	8	&	9	&	10	&	11	&	12	&	13	&	14	&	15	& 16     & $n^{med}_g$	\\\hline
1	&	1	&		&		&		&		&		&		&		&		&		&		&		&		&		&	&      1	\\
2	&	1	&	1	&		&		&		&		&		&		&		&		&		&		&		&		&       & 	   2	\\
3	&	1	&	1	&	1	&		&		&		&		&		&		&		&		&		&		&		&	&      3	\\
4	&	1	&	1	&	1	&	1	&		&		&		&		&		&		&		&		&		&		& 	&      4	\\
5	&	1	&	2	&	2	&	1	&	1	&		&		&		&		&		&		&		&		&		&       &       	7	\\
6	&	1	&	2	&	3	&	2	&	1	&	1	&		&		&		&		&		&		&		&		&	&        10	\\
7	&	1	&	2	&	4	&	2	&	2	&	1	&	1	&		&		&		&		&		&		&		&       &    	13	\\
8	&	1	&	3	&	5	&	4	&	4	&	2	&	1	&	1	&		&		&		&		&		&		&	&       21	\\
9	&	1	&	3	&	7	&	5	&	6	&	4	&	2	&	1	&	1	&		&		&		&		&		&	&      30	\\
10	&	1	&	3	&	8	&	8	&	9	&	4	&	4	&	2	&	1	&	1	&		&		&		&		&	&      41	\\
11	&	1	&	4	&	10	&	10	&	14	&	7	&	7	&	4	&	2	&	1	&	1	&		&		&		&	&       61	\\
12	&	1	&	4	&	12	&	13	&	19	&	12	&	10	&	7	&	4	&	2	&	1	&	1	&		&		&	&       86	\\
13	&	1	&	4	&	14	&	16	&	25	&	18	&	17	&	9	&	7	&	4	&	2	&	1	&	1	&		&	&      119	\\
14	&	1	&	5	&	16	&	22	&	35	&	25	&	26	&	16	&	12	&	7	&	4	&	2	&	1	&	1	&	&      173	\\
15	&	1	&	5	&	19	&	24	&	45	&	37	&	39	&	24	&	47	&	27	&	15	&	4	&	2	&	1	&   1	&     291	\\\hline
\end{tabular}
\end{center}}
\caption{Number of MED-semigroups with given gender $g$ and multiplicity $m$.\label{med0to15}}
\end{table}
Analogous results may be stated for symmetric numerical semigroups or maximal embedding dimension numerical semigroups, using the correspondence between the corresponding family and the polytope describing the \'Apery set.

In the following section we exploit the use of short generating functions for the smallest multiplicities, where the polytopes are more manageable and formulas for the number of integer points inside them can be obtained for any genus or Frobenius number by using differente techniques.
\section{Numerical Semigroups with multipicities three and four}
\label{sec:cases}
In this section we analyze, using the above methodology, those numerical
semigroups with multiplicities three and four. For these families we obtain the
short generating functions defining the polytopes characterizing these semigroups with arbitrary genus or Frobenius number. Taking advantage of the 2-dimensional and 3-dimensional geometry of the corresponding
 bodies, we give formulas for the number of these numerical semigroups.
\subsection{Multiplicity Three}
Let $\mathcal{S}_{3,g}$ be the set of numerical semigroups with multiplicity $m=3$ and genus $g$. The following result depicts the number of numerical semigropus in $\mathcal{S}_{3,g}$.  It is based on the structure of the 2-dimensional polytope that characterizes the set $\mathcal{S}_{3,g}$.
\begin{theo}
The number of numerical semigroups of multiplicity $3$ and genus $g$ is $\left\lceil \dfrac{g+1}{3} \right\rceil$.
\end{theo}
\begin{proof}
Let $P_{3,g}$ be the polytope \eqref{eq:polytope} when $m=3$, that is, $P_{3,g} = \{(x,y): 2x-y\ge 0, -x+2y\ge -1, x+y=g, x\ge 1, y \ge 1\}$.

Since the integer points inside $P_{3,g}$ must be in the line $x+y=g$ we can reduce this polytope to the 1-dimensional polytope $P_{3,g}^\prime = \{x: \frac{g}{3} \le x \le \frac{2g+1}{3}\}$. Then, since $y = g-x$ is
integer if and only if $x$ is integer, the number of integer points inside $P_{3,g}$ is $\# (P^\prime_{3,g} \cap \Z) = \left\lfloor \frac{2g+1}{3} \right\rfloor - \left\lceil \frac{g}{3} \right\rceil = \left\lceil \frac{g+1}{3} \right\rceil$.
\end{proof}
The above result coincides with Corollary 10 in \cite{rosales05} that was proved following a different strategy.
\begin{remark}[Maximal embedding dimension numerical semigroups with multiplicity $3$]
An analogous analysis can be done for MED-semigroups using the polytope characterization for these semigroups \eqref{eq:meds}. In this case, the number of MED-semigroups of multiplicity $3$ and any genus $g\ge 2$ is $\left\lceil \dfrac{g-1}{3} \right\rceil$.

Let $\mathcal{S}_{3,g}$ be the set of numerical semigroups of genus
$g$ and multiplicity $3$ and $\mathcal{S}^{me}_{3,g}$ the set of
MED-semigroups of genus $g$ and multiplicity $3$.
\begin{enumerate}
\item The set $\mathcal{S}_{3,g}$ coincides with the set $\mathcal{S}^{me}_{3,g}$ if and only if $g \equiv 2 \pmod 3$.
\item If $g \not\equiv 2 \pmod 3$, then, $\# \mathcal{S}^{me}_{3,g} = \mathcal{S}^{me}_{3,g} + 1$, i.e., there exists only
 one numerical semigroup of genus $g$ and multiplicity $3$ that does not have maximal embedding dimension.
\end{enumerate}
\end{remark}

Note that the above remark is not surprising. If $S$ is an embedding dimension three numerical semigroup that is not of maximal embedding dimension (equals to $3$ in this case) and genus $g$, then $S=\langle 3, n \rangle$ for some $n>3$. Since in this setting $\g(S)=\frac{\Frob(S)+1}{2} = \frac{1}{2}\,(3-1)(n-1) = n-1$, we get that $S=\langle 3, g+1\rangle$. If $g\equiv 2 \pmod 3$, $g+1$ is a multiple of $3$, and $S$ is not a numerical semigroup, and then there are no numerical semigroups with multiplicity $3$ that are not of maximal embedding dimension. Otherwise ($g\equiv 0 \text{ or } 1 \pmod 3$) there is a unique numerical semigroup that is not of maximal embedding dimension.

In the following we count numerical semigroups with multiplicity $3$ and arbitrary Frobenius number. In this case, the set of numerical semigroups of multiplicity $3$ and Frobenius number $F$ is described by the system of diophantine inequalities \eqref{eq:polyFrob} when $m=3$.

We distinguish two cases, attending to the two possible choices for $F$: $F\equiv 1 \pmod 3$ and $F \equiv 2 \pmod 3$.
\begin{itemize}
\item If $F\equiv 1 \pmod 3$, then $k^* = 1$, and the
System \eqref{eq:polyFrob} is reduced by substituting $x_1$ by $\frac{F+2}{3}$, the constraints for $x_2$ are $x_2\leq \frac{F+1}{3}$ and $x_2 \geq \frac{F-1}{6}$. Then, the number of possibilities for $x_2$ are $\lfloor\frac{F+1}{3}\rfloor - \lceil \frac{F-1}{6} \rceil + 1$.
\item If $F\equiv 2\pmod 3$, then $k^* = 2$, after substituting $x_2$ by $\frac{F+1}{3}$ the constraints for $x_1$ become $x_1\leq \frac{F+2}{3}$ and $x_1 \geq \frac{F+1}{6}$. Then, the number of possibilities for $x_1$ are $\lfloor\frac{F+2}{3}\rfloor - \lceil \frac{F+1}{6} \rceil + 1 = \lfloor\frac{F+1}{3}\rfloor - \lceil \frac{F-1}{6} \rceil + 1$, coinciding with the case when $F \equiv 1 \pmod 3$.
\end{itemize}
\begin{proposition}
The number of numerical semigroups with multiplicity $3$ and Frobenius number $F$ (that is not a multiple of $3$) is $\lfloor\frac{F+1}{3}\rfloor - \lceil \frac{F-1}{6} \rceil + 1$.
\end{proposition}
By using the same methodology, but applied to MED semigroups, we obtain the following result.
\begin{proposition}
The number of MED semigroups with multiplicity $3$ and Frobenius number $F$ (that is not a multiple of $3$) is $\lfloor\frac{F+1}{3}\rfloor - \lceil \frac{F+2}{6} \rceil + 1$.
\end{proposition}
\begin{remark}[Fixing genus and Frobenius number]
\label{fixgF:m3}
With this procedure, adding to the polytope describing numerical semigroups with fixed genus, $g$, those constraints that fix the Frobenius number, $F$,
 we count the number of numerical semigroups with both values fixed. In the case when the multiplicty is $3$, we obtain that all these polytopes (varying $F$ and $g$) have a unique integer
point inside them. This fact reproves the well-known result that states that a numerical semigroup with multiplicity $3$ is completely determined by its genus and its Frobenius number {\rm(see \cite{rosales05})}. Actually, the unique lattice point inside that polytope is $(\frac{F+2}{3}, \frac{3g-F-2}{3})$ if $F\equiv 1 \pmod 3$ or $( \frac{3g-F-1}{3}, \frac{F+1}{3})$ if $F\equiv 2 \pmod 3$ being the unique coordinates of the elements in the Ap\'ery set $\Ap(S, m)=\{0, F+3, 3g-F\}$. Then, the unique numerical semigroup with multiplicity $3$, genus $g$ and Frobenius number $F$ is $S = \langle \Ap(S, m), m\rangle = \langle 3, F+3, 3g-F \rangle$.
\end{remark}
\subsection{Multiplicity Four}
Let $\mathcal{S}_{4,g}$ be the set of numerical semigroups with multiplicity $m=4$ and genus $g$. In the analysis for numerical semigroups with multiplicity $3$, we reduce the problem to dimension $1$ since the polytope characterizing these semigroups is a $2$-dimensional polyhedron and the equality constraint that fix the genus allows us to change one of the variables by a linear expression in the other.
In this subsection we study, arguing analogously, the multiplicity four case. The main difference between this setting and the previous one is that in this case ($m=4$) we are able to reduce the problem to a $2$-dimensional polytope and generating functions tools are needed to count the integer points inside it. This is also the reason why to analyze larger multiplicities become much more complicated in $3$ or more dimensions (see comment at the end of this section).

We distinguish different ranges of genus, since the geometry of the polytopes characterizing the numerical semigroups are slightly different.
\begin{proposition}
Let $P_{4,g} = \{(x,y, z): 2x-y\ge 0, x+y-z \ge 0, -x+y+z \ge -1, -y + 2z \ge -1, x+y+z=g, x\ge 1, y \ge 1, z\ge 1\}$.
\begin{itemize}
 \item If $g \in [3, 8]$ then $\# (P_{4,g} \cap \Z^3) = \left\{\begin{array}{rl} 1 & \mbox{if $g=3$},\\
							  3 & \mbox{if $g=4$},\\
							  4 & \mbox{if $g=5$},\\
							  6 & \mbox{if $g=6$},\\
							  7 & \mbox{if $g=7$},\\
						          9 & \mbox{if $g=8$}.
                                                        \end{array}\right.$
 \item If $g\ge 9$, then $\# (P_{4,g} \cap \Z^3)= \# (T_A(g) \cap \Z^2) + \# (T_B(g) \cap \Z^2) + \# (R(g) \cap \Z^2) - \# (T_C(g) \cap \Z^2)$
where
\begin{itemize}
\item $T_A(g) = \{ (x,y) \in \R^2: 3x+y \ge g, x\le \frac{2g+1}{8}, y \le \frac{g}{2}\}$,
\item $T_B(g) = \{ (x,y) \in \R^2: x+3y \ge g-1, x\le \frac{g+1}{2}, y \le \frac{2g-3}{8}\}$,
\item $R(g) = \{ (x,y) \in \R^2: \frac{2g+1}{8} \le x\le \frac{g+1}{2}, \frac{2g-3}{8} \le y \le \frac{g}{2}\}$, and
\item $T_C(g) = \{ (x,y) \in \R^2: x+y \ge g, x \le \frac{g+1}{2}, y \le \frac{g}{2}\}$.
\end{itemize}
\end{itemize}
\end{proposition}
\begin{proof}
 The polytope $P_{4,g}$ in $\R^3$ can be projected in a lower dimension space by using the constraint $x+y+z = g$. By substituting each $y$ by the expression $g-x-z$, we obtaing the following $2$-dimensional polytope:
$P^\prime = \{(x,y): 3x+y\ge g, z \le \frac{g}{2}, x \le \frac{g+1}{2}, x+3y \ge g-1, x+y \le g-1, x\ge 1, y \ge 1\}$.
\begin{itemize}
\item Clearly, for $g=3$, the number of integer points inside $P_{4,3}$ is $1$.
 \item For $g=4$, by using \texttt{barvinok} the generating function
   for the corresponding polytope is
$$
\dfrac{z_2^{4}}{(1-{\frac {z_{{2}}}{z_{{1}}}})(1-z_{{2}})} -\dfrac{z_1^{3}z_2}{(1-\frac{z_2}{z_1})( 1-z_1)} + \frac{z_1z_2}{(1-z_2)(1-z_1)}.
$$
By substituting in the above expression $z_1=\xi$ and $z_2=\xi^2$, and taking limit when $\xi \rightarrow 1$ we obtain that the number of integer points inside that polytope is $3$.
 \item For $g=5$, \texttt{barvinok} output is 
$$
 \dfrac{z_1^{2} z_2^{3} - z_1^{4} z_2}{( 1-\frac{z_2}{z_1})( 1-z_1)}+ \dfrac{z_1z_2^2-z_2^5}{(1-\frac{z_2^3}{z_1})(1-z_2)}+\dfrac{z_1^2z_2-z_1z_2^{3}}{(1-z_2)(1-z_1)}.
$$
Proceeding as above we get that the number of integer points inside that polytope is $4$.
 \item For $g=6$, invoking \texttt{barvinok} again, the generating
   function for the corresponding polytope is
$$
\frac{-\frac{z_2^2}{z_1}+z_1^{2}z_2}{( 1-\frac{z_2}{z_1^{3}})( 1-z_1)}+  \frac{-z_1z_2^{2}-z_1^{3}z_2^{3}}{(1-\frac {z_2}{z_1})(1-z_2)} +
\frac{z_1z_2^{2}+z_1^{2}z_2^{4}}{( 1-\frac{z_2}{z_1})( 1-z_1)}+ \frac{-z_2^{6}+z_1z_2^{3}}{( 1-\frac{z_2^{3}}{z_1})( 1-z_2)}+\frac{-z_1z_2^{4}-z_1^4z_2}{
(1-z_{{2}})( 1-z_{{1}})}.
$$
Repeating the above limit procedure, the number of integer points inside that polytope is $6$.
 \item For $g=7$, the generating function is
$$
\frac{-{z_{{2}}}^{2}+{z_{{1}}}^{3}z_{{2}}}{( 1-{\frac {z_{{2}}}{{z_{{1}}}^{3}}})( 1-z_{{1}})}+ \frac{-z_{{1}}{z_{{2}}}^{3}-{z_{{1}}}^{4}{z_{{2}}}^{3}}
{( 1-{\frac {z_{{2}}}{z_{{1}}}})( 1-z_{{2}})} + \frac{{z_{{1}}}^{2}{z_{{2}}}^{2}+{z_{{1}}}^{3}{z_{{2}}}^{4}}{( 1-{\frac {z_{{2}}}{z_{{1}}}})( 1-z_{{1}})}
+{\frac {-{z_{{1}}}^{2}{z_{{2}}}^{4}-{z_{{1}}}^{5}z_{{2}}}{ \left( 1-z_{{2}} \right)  \left( 1-z_{{1}} \right)}}.
$$
Being $7$ the number of integer points inside that polytope.
 \item For $g=8$, \texttt{barvinok} yields 
$$\frac{-z_{{1}}{z_{{2}}}^{2}+{z_{{1}}}^{4}z_{{2}}}{( 1-{\frac {z_{{2}}}{{z_{{1}}}^{3}}})( 1-z_{{1}})}+ \frac{-{z_{{1}}}^{2}{z_{{2}}}^{2}-{z_{{1}}}^{4}{z_{{2}}}^{4}}
{( 1-{\frac{z_{{2}}}{z_{{1}}}})( 1-z_{{2}})}+ \frac{z_1^{2}z_2^{2}+z_1^{3}z_2^5}{( 1-{\frac {z_{{2}}}{z_{{1}}}})( 1-z_{{1}})}
+\frac{-z_{{1}}{z_{{2}}}^{5}+{z_{{1}}}^{2}{z_{{2}}}^{2}}{( 1-{\frac {{z_{{2}}}^{3}}{z_{{1}}}})( 1-z_{{2}})}
+{\frac {-{z_{{1}}}^{2}{z_{{2}}}^{5}-{z_{{1}}}^{5}z_{{2}}}{ \left( 1-z_{{2}} \right)  \left( 1-z_{{1}} \right) }}.
$$
And then, taking limit, as above, we obtain that the number of integer points inside that polytope is $9$.
 \item If $g\ge 9$, the polytope can be decomposed as the disjoint union $P^\prime \cap \Z^2= (T_A(g) \cap \Z^2) \cup (T_B(g) \cap \Z^2) \cup
 (R(g)\backslash T_C(g) \cap \Z^2)$. Let $(x,y) \in P^\prime \cap \Z^2$, if $x \le \frac{2g+1}{8}$, then, because of the constraint $3x+y \ge g$, $(x,y)$
is in $T_A(g) \cap \Z^2$. If $y\le \frac{2g-3}{8}$, then, in view of $x+3y \ge g-1$, $(x,y)$ is in $T_B(g) \cap \Z^2$. Otherwise,
it is clear that the point $(x,y)$ is in $R(g) \backslash T_C(g)$. Note that since the vertex $(\frac{2g+1}{8}, \frac{2g-3}{8})$ is not an integer vector there are no integer points in any pairwise
 intersections of $T_A(g), T_B(g)$ and $R(g)\backslash T_C(g)$. Figure \ref{fig:polym4} shows this decomposition.
\end{itemize}
\begin{figure}[H]
 \begin{center}
  \includegraphics[scale=0.7]{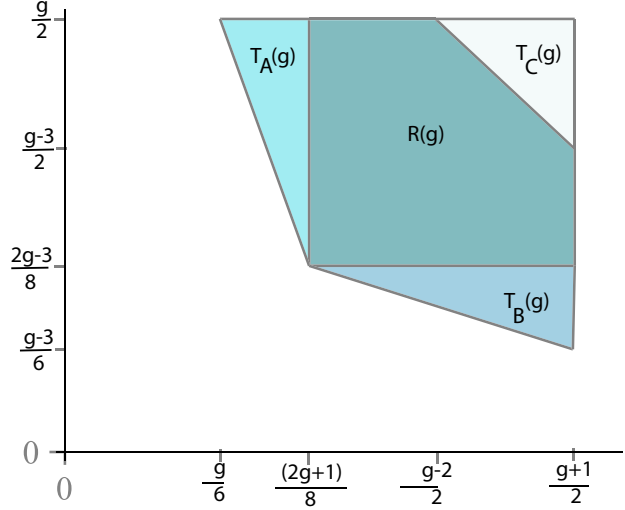}
 \end{center}
\caption{Decomposition of the polytopes for $g\ge 9$.\label{fig:polym4}}
   \vspace*{-.5cm}
\end{figure}
\end{proof}
In view of the above result, we can count the integer points inside
each of the parts (triangles and rectangle) desctibing $P_{4,g}$ with
$g\ge 9$, and adding all of them to obtain the number of points inside
the whole polytope. By using the software \texttt{barvinok} for
parametric polytopes (with parameter $g\ge 9$), we obtain the
following formulas for the number of points inside each of the
triangles. 

$
\# (T_A(g) \cap \Z^2)  = \frac{1}{2}\left(\left\lfloor \frac{g+5}{6} \right\rfloor - \left\lfloor\frac{g}{4} \right\rfloor -1\right) \left(2g-3 \left\lfloor \frac{g}{4} \right\rfloor -2 \left\lfloor \frac{g}{2} \right\rfloor - 2 - 3 \left\lfloor \frac{g+5}{6} \right\rfloor \right),
$

$
\# (T_B(g) \cap \Z^2)   = \frac{1}{2}\left(\left\lfloor \frac{g+2}{6} \right\rfloor - \left\lfloor\frac{g+2}{4} \right\rfloor \right) \left(-3 \left\lfloor \frac{g+2}{4} \right\rfloor +2 \left\lfloor \frac{g}{2} \right\rfloor- 3 \left\lfloor \frac{g+2}{6} \right\rfloor -1\right),
$

$
\# (T_C(g) \cap \Z^2)  = 1 \mbox{ (this is easily checkable analyzing the vertices of the rectangle and isosceles triangle $T_C(g)$).}
$

The generating function for a box $H=[a,b] \times [c,d] \subseteq
\R_+^2$ is
$$
f(H; z_1, z_2) = \dfrac{z_1^{\lceil a \rceil} z_2^{\lceil c \rceil}}{(1-z_1)(1-z_2)} + \dfrac{z_1^{\lfloor b \rfloor} z_2^{\lceil a \rceil}}{(1-z_1^{-1})(1-z_2)} + \dfrac{z_1^{\lceil a \rceil} z_2^{\lfloor d \rfloor}}{(1-z_1)(1-z_2^{-1})} + \dfrac{z_1^{\lfloor b \rfloor} z_2^{\lfloor d \rfloor}}{(1-z_1^{-1})(1-z_2^{-1})}.
$$
Set $z_1=\xi$ and $z_2=\xi^2$, and take limit when $\xi \rightarrow 1$. We obtain the well-known formula for the number of integer points inside the box $H$:
$$
(1 + \lfloor b \rfloor + \lceil a \rceil)(1 + \lfloor d \rfloor + \lceil b \rceil).
$$
By applying this formula to the rectangle $R(g)$ we obtain the
following descriptions for the numbers of points inside them.
$$
\# (R(g) \cap \Z^2) = \left(\left\lfloor \frac{g+1}{2} \right\rfloor - \left\lceil \frac{2g+1}{8}\right\rceil +1\right)\left(\left\lfloor \frac{g}{2} \right\rfloor - \left\lceil \frac{2g-3}{8} \right\rceil + 1\right).
$$
The following result summarizes the above reasonings.
\begin{theo}
 Let $\mathcal{S}_{4,g}$ be the set of numerical semigroups of
 multiplicity $4$ and genus $g \ge 9$. Then the cardinality of
 $\mathcal{S}_{4,g}$ is
$-g + \frac{5}{2} \lfloor \frac{g}{4} \rfloor + +\lfloor \frac{g}{2} \rfloor+\frac{1}{2}\lfloor \frac{g+5}{6} \rfloor+\lfloor \frac{g+5}{6} \rfloor g -\lfloor \frac{g+5}{6} \rfloor\lfloor \frac{g}{2} \rfloor-\frac{3}{2}\lfloor \frac{g+5}{6} \rfloor^2-\lfloor \frac{g}{4} \rfloor g+\frac{3}{2} \lfloor \frac{g}{4} \rfloor^2+\lfloor \frac{g}{4} \rfloor\lfloor \frac{g}{2} \rfloor-\frac{1}{2}\lfloor \frac{g+2}{6} \rfloor+\lfloor \frac{g+2}{6} \rfloor\lfloor \frac{g}{2} \rfloor-\frac{3}{2}\lfloor \frac{g+2}{6} \rfloor^2+\frac{1}{2} \lfloor \frac{g+2}{4} \rfloor+\frac{3}{2} \lfloor \frac{g+2}{4} \rfloor^2-\lfloor \frac{g+2}{4} \rfloor\lfloor \frac{g}{2} \rfloor+\lfloor \frac{g+1}{2} \rfloor\lfloor \frac{g}{2} \rfloor-\lfloor \frac{g+1}{2} \rfloor\lceil \frac{2g-3}{8} \rceil+\lfloor \frac{g+1}{2} \rfloor-\lceil \frac{2g-7}{8} \rceil \lfloor \frac{g}{2} \rfloor+\lceil \frac{2g-7}{8} \rceil\lceil \frac{2g-3}{8} \rceil-\lceil \frac{2g-7}{8} \rceil$.
\end{theo}
If we check the formula in the above result for those genus less than $9$, we see that it is also verified for those genus $g \in \{4,5,6,7,8\}$. For $g=3$ is is known that there exists only one numerical semigroup with multiplicity $4$ and genus $3$.
\begin{corollary}
 Let $\mathcal{S}_{4,g}$ be the set of numerical semigroups of multiplicity $4$ and genus $g$. If $g=3$, then  $\mathcal{S}_{4,g} =1$, otherwise  the cardinality of $\mathcal{S}_{4,g}$ is $-g + \frac{5}{2} \lfloor \frac{g}{4} \rfloor + +\lfloor \frac{g}{2} \rfloor+\frac{1}{2}\lfloor \frac{g+5}{6} \rfloor+\lfloor \frac{g+5}{6} \rfloor g -\lfloor \frac{g+5}{6} \rfloor\lfloor \frac{g}{2} \rfloor-\frac{3}{2}\lfloor \frac{g+5}{6} \rfloor^2-\lfloor \frac{g}{4} \rfloor g+\frac{3}{2} \lfloor \frac{g}{4} \rfloor^2+\lfloor \frac{g}{4} \rfloor\lfloor \frac{g}{2} \rfloor-\frac{1}{2}\lfloor \frac{g+2}{6} \rfloor+\lfloor \frac{g+2}{6} \rfloor\lfloor \frac{g}{2} \rfloor-\frac{3}{2}\lfloor \frac{g+2}{6} \rfloor^2+\frac{1}{2} \lfloor \frac{g+2}{4} \rfloor+\frac{3}{2} \lfloor \frac{g+2}{4} \rfloor^2-\lfloor \frac{g+2}{4} \rfloor\lfloor \frac{g}{2} \rfloor+\lfloor \frac{g+1}{2} \rfloor\lfloor \frac{g}{2} \rfloor-\lfloor \frac{g+1}{2} \rfloor\lceil \frac{2g-3}{8} \rceil+\lfloor \frac{g+1}{2} \rfloor-\lceil \frac{2g-7}{8} \rceil \lfloor \frac{g}{2} \rfloor+\lceil \frac{2g-7}{8} \rceil\lceil \frac{2g-3}{8} \rceil-\lceil \frac{2g-7}{8} \rceil$.
\end{corollary}
In the following we describe the case when the multiplicity is four and the Frobenius number, $F$, varies.

The set of numerical semigroups of multiplicity $4$ and Frobenius number $F$ is described by the of diophantine inequalities in \eqref{eq:polyFrob} when $m=4$. We distinguish three cases. 
\begin{itemize}
\item If $F\equiv 1 \pmod 4$, then $k^* = 1$, and the
System \eqref{eq:polyFrob} is reduced, after substituting $x_1$ by $\frac{F+3}{4}$:
\begin{eqnarray}
x_2, x_3  \geqslant&1, \nonumber\\
4\,x_3-4\,x_2 \leqslant& F+3, \nonumber\\
4\,x_2+4\,x_3 \geqslant&F-1,\nonumber\\
-x_2+2x_3 \geqslant&-1,\nonumber\\
4\,x_2 &\leq F + 2,\nonumber\\
4\,x_3 &\leq F + 1,\nonumber\\
x_2, x_3\in \Z. & \nonumber
\end{eqnarray}
Using \texttt{barvinok} the number of solutions for the above
parametric polytope is
{\small
$$
\left\{ \begin{array}{ll}
\left(\frac{F-1}{4}\right)^2 & \mbox{if $5 \leq F \leq 9$},\\
\frac{F^2-14F+141}{16} & \mbox{if $13 \leq F \leq 17$},\\
\frac{-3}{2}\left\lfloor \frac{F+1}{12} \right\rfloor^2  + \frac{1}{4}\left\lfloor \frac{F+1}{12} \right\rfloor\,F \frac{-3}{4}\left\lfloor \frac{F+1}{12} \right\rfloor + \left\lfloor \frac{F+5}{8} \right\rfloor^2 -\frac{1}{4}\left\lfloor \frac{F+5}{8} \right\rfloor\,F +\frac{1}{4}\left\lfloor \frac{F+5}{8} \right\rfloor + \frac{5F}{6} - \frac{11}{32} + \frac{F^2}{32}& \mbox{if $F \geq 21$}.\end{array}\right.
$$}
\item If $F\equiv 2 \pmod 4$, then $k^* = 2$, and the
System \eqref{eq:polyFrob} is reduced, after substituting $x_2$ by $\frac{F+2}{4}$:
\begin{eqnarray}
x_1, x_3  \geqslant&1, \nonumber\\
4\,x_3-4\,x_1 \leqslant& F+2, \nonumber\\
4\,x_3+4\,x_1 \geqslant&F-2,\nonumber\\
F- 2 \leqslant 8\,x_3 \leqslant&F+2,\nonumber\\
F+2 \leqslant 8\,x_1 &\leqslant 2F + 6,\nonumber\\
x_1, x_3\in \Z. & \nonumber
\end{eqnarray}
\texttt{barvinok} output for the number of solutions for the above
parametric polytope is
$$
\left(\frac{F}{4}-\left\lfloor \frac{F+1}{8}\right\rfloor + \frac{1}{2}\right)\,\left(\frac{F}{4}-\left\lfloor \frac{F+5}{8}\right\rfloor + \frac{1}{2}\right).
$$
\item If $F\equiv 3 \pmod 4$, then $k^* = 3$, and the
System \eqref{eq:polyFrob} is reduced, after substituting $x_3$ by $\frac{F+1}{4}$:
\begin{eqnarray}
x_1, x_2  \geqslant&1, \nonumber\\
2\,x_1 - x_2 \geqslant& 0, \nonumber\\
4\,x_1+4\,x_2 \geqslant& F+1, \nonumber\\
4\,x_2-4\,x_1 \leqslant&F+5,\nonumber\\
4\,x_1 &\leq F + 3,\nonumber\\
4\,x_2 &\leq F + 2,\nonumber\\
x_1, x_2\in \Z. & \nonumber
\end{eqnarray}
Invoking again \texttt{barvinok}, the number of solutions for the
above parametric polytope is
$$
\frac{1}{32} F^2 + \frac{7}{16} F -\frac{19}{32} + \frac{3}{4} \left\lfloor \frac{F+1}{8} \right\rfloor + \left \lfloor \frac{F+1}{8} \right\rfloor^2 + \frac{1}{4} \left \lfloor \frac{F}{12} \right\rfloor F - \frac{9}{4} \left \lfloor \frac{F}{12} \right\rfloor - \frac{3}{2} \left \lfloor \frac{F}{12} \right\rfloor^2 - \frac{1}{4} \left \lfloor \frac{F+1}{8} \right\rfloor F.
$$
\end{itemize}
\begin{theo}
The number of numerical semigroups with multiplicity $4$ and
Frobenious number $F$ is
{\scriptsize
$$
\left\{\begin{array}{rl}\left(\frac{F-1}{4}\right)^2 & \mbox{if $F\equiv 1 \pmod 4$ and $5 \leq F \leq 9$},\\
\frac{F^2-14F+141}{16} & \mbox{if $F\equiv 1 \pmod 4$ and $13 \leq F \leq 17$},\\
\frac{-3}{2}\left\lfloor \frac{F+1}{12} \right\rfloor^2  + \frac{1}{4}\left\lfloor \frac{F+1}{12} \right\rfloor\,F -\frac{3}{4}\left\lfloor \frac{F+1}{12} \right\rfloor + \left\lfloor \frac{F+5}{8} \right\rfloor^2 -\frac{1}{4}\left\lfloor \frac{F+5}{8} \right\rfloor\,F +\frac{1}{4}\left\lfloor \frac{F+5}{8} \right\rfloor + \frac{5F}{6} - \frac{11}{32} + \frac{F^2}{32}& \mbox{if $F\equiv 1 \pmod 4$ and $F \geq 21$},\\
\left(\frac{F}{4}-\left\lfloor \frac{F+1}{8}\right\rfloor + \frac{1}{2}\right)\,\left(\frac{F}{4}-\left\lfloor \frac{F+5}{8}\right\rfloor + \frac{1}{2}\right) & \mbox{if $F\equiv 2 \pmod 4$},\\
\frac{F^2}{32} + \frac{7F}{16} -\frac{19}{32} + \frac{3}{4} \left\lfloor \frac{F+1}{8} \right\rfloor + \left \lfloor \frac{F+1}{8} \right\rfloor^2 + \frac{1}{4} \left \lfloor \frac{F}{12} \right\rfloor F - \frac{9}{4} \left \lfloor \frac{F}{12} \right\rfloor - \frac{3}{2} \left \lfloor \frac{F}{12} \right\rfloor^2 - \frac{1}{4} \left \lfloor \frac{F+1}{8} \right\rfloor F
& \mbox{if $F\equiv 3 \pmod 4$},\\
0 & \mbox{otherwise}.
\end{array}\right.
$$}
\end{theo}
\begin{remark}[Fixing genus and Frobenius number]
It is well-known that fixing genus and Frobenius number is not enough to determine, in general, a
numerical semigroup with multiplicity $4$ (see Remark \ref{fixgF:m3} for the analysis of the case of multiplicity $3$). If we add to the polytope encoding the numerical semigroups with given genus, $g$, those constraints fixing
 the Frobenius number, $F$, we obtain this output in \verb"barvinok":
$$
\left\{\begin{array}{rl} \frac{F+3}{2} - \left\lfloor\frac{2g+F+5}{6}\right\rfloor & \mbox{if $F\equiv 1 \pmod 4$, $5F-8g\ge 5$, and $2g-F\ge 5$},\medskip\\
g - \left\lfloor\frac{2g+F+5}{6}\right\rfloor & \mbox{if $F\equiv 1 \pmod 4$, $4g-F\ge 23$, $2g-F \ge 1$ and $2g-F\le 3$},\medskip\\
\frac{3F-4g+5}{4} & \mbox{if $F\equiv 1 \pmod 4$, $5F-8g\le 1$, $4g-3F\le 1$ and $2g-F\ge 5$},\\
\frac{4g-F-7}{4} & \mbox{if $F\equiv 1 \pmod 4$, $4g-F\ge 11$, $4g-F\le 19$, $2g-F\ge 1$and $2g-F\le 3$},\\
\frac{2g-F}{2} & \mbox{if $F\equiv 2 \pmod 4$, $8g-5F\le 2$, $F\ge 14$, and $2g-F\ge 2$},\\
\frac{3F-4g+6}{4}  & \mbox{if $F\equiv 2 \pmod 4$, $8g-5F\ge 6$, $4g-3F\le 2$ and $2g-F\ge 6$},\\
\frac{F-2}{4} & \mbox{if $F\equiv 2 \pmod 4$, $8g-5F\ge 6$, $F\ge 6$, and $2g-F\le 4$},\medskip\\
\frac{4g-F-2}{4} - \left\lfloor\frac{F+6}{8}\right\rfloor & \mbox{if $F\equiv 2 \pmod 4$, $8g-5F\le 2$, $F\le 10$, and $2g-F\ge 2$,}\medskip\\
\frac{3F-4g+7}{4}& \mbox{if $F\equiv 3 \pmod 4$, $8g-5F\ge 9$, $4g-3F\le 3$ and $2g-F\ge 5$,}\medskip\\
\frac{F+3}{2} - \left\lfloor\frac{2g+F+5}{6}\right\rfloor & \mbox{if $F\equiv 3 \pmod 4$, $8g-5F\le 5$, and $2g-F\ge 5$,}\medskip\\
g - \left\lfloor\frac{2g+F+5}{6}\right\rfloor & \mbox{if $F\equiv 3 \pmod 4$, $8g-5F\le 5$, $4g-F\ge 9$, $2g-F\ge 1$, and $2g-F\le 3$}\\
0 & \mbox{otherwise}.
       \end{array}\right.
$$
Then, it is clear that different multiplicity four numerical semigroups can have the same Frobenius number and genus. If we analyze the three cases attending to the different possibilities for the embedding dimension of the numerical semigroups we see that this result is not unusual. Let $S$ be a numerical semigroup with multiplicity $4$, genus $g$, Frobenius number $F$, and embedding dimension ${\rm e}(S)$.
\begin{enumerate}
 \item If ${\rm e}(S)=2$, then $S=\langle 4 < n_2 \rangle$. From Selmer's formulas, $n_2=\frac{F+4}{3}$ and $g=\frac{F+1}{4}$, and then, $S$ is completely determined by $F$, being $S=\langle 4, \frac{F+4}{3}\rangle$.
 \item If ${\rm e}(S)=3$, then $S=\langle 4 < n_2 < n_3 \rangle$. This case, is particulary difficult and is deeply studied in \cite{rosales09}, where the authors prove that $S$ is completely characterized by $F$ and $g$.
 \item If ${\rm e}(S)=4$ (maximal embedding dimension case), then
   $S=\langle 4 < n_2 < n_3 < n_4 \rangle$. From Selmer's formulas,
   $n_4=F+4$ and $g=\frac{n_2 + n_3 + F + 4}{4} - \frac{4 - 1}{2}$,
   being $n_2$ and $n_3$ all the integer solutions for $n_2 + n_3 =
   4g-F+2$ with $4 < n_2 < n_3 < F+4$. This diophantine equation has,
   in general, more than one solution, actually, invoking
   \verb"barvinok" for counting these solutions we obtain that the
   number of solutions is
     {\small
$$ 
\left\{\begin{array}{rl}
2g-\left\lfloor \frac{F}{2} \right\rfloor - 4 & \mbox{if $2g-F \le 2$ and $4g-F \ge 9$,}\\
2F-2g-\left\lfloor \frac{F}{2} \right\rfloor + 2 & \mbox{if $2g-F \ge 3$ and $4g-3F \le 3$.}\\
\end{array}\right.
$$}
However, in this system is not taken into account that
$\gcd(n_i, n_j)=1$ for $i, j=1, \ldots, 4$, $i\neq j$, and then, the
above formula is only a upper bound for the number of MED-semigroups
with multiplicity $4$, genus $g$ and Frobenius number $F$. In case we
want the exact number of these semigroups, we need to use the
corresponding formula proceeding as in the general case, where
\verb"barvinok" output is

$
\left\{\begin{array}{rl}
 \frac{F+1}{2} - \left\lfloor \frac{2g+F+1}{6}  \right\rfloor & \mbox{if $F\equiv 1 \pmod 4$ and  $-8g+5F\ge 1$ and $2g-F \le 3$}\\
\frac{3F+5}{4} - g & \mbox{if $F\equiv 1 \pmod 4$ and  $8g-5F\ge 3$ and $4g-F \ge 15$ and $4g-3F\le 1$,}\\
 \frac{F}{2} - \left\lfloor \frac{2g+F}{6}  \right\rfloor & \mbox{if $F\equiv 2 \pmod 4$ and  $8g-5F\le 6$ and $2g-F \ge 2$ and $4g-F\ge 14$,}\\
\frac{3F+6}{4} - g & \mbox{if $F\equiv 2 \pmod 4$ and  $8g-5F\ge 10$ and $4g-3F\le 2$,}\\
 \frac{F+3}{2} - \left\lfloor \frac{2g+F+5}{6}  \right\rfloor & \mbox{if $F\equiv 3 \pmod 4$ and  $8g-5F\le 1$ and $2g-F \ge 5$,}\\
g - \left\lfloor \frac{2g+F+5}{6} \right\rfloor & \mbox{if $F\equiv 3 \pmod 4$ and  $8g-5F\le 1$ and $4g-F\ge 13$ and $1 \le 2g-F\le 3$,}\\
\frac{3F+7}{4} - g & \mbox{if $F\equiv 3 \pmod 4$ and  $8g-5F\ge 5$ and $4g-3F\le 3$ and $2g-F\ge 5$,}\\
2 & \mbox{if $F=5$ and $g=7$,}\\
0 & \mbox{otherwise.}
\end{array}\right.
$
\end{enumerate}
\end{remark}
The methodology applied in this paper for numerical semigroups with multiplicities $3$ and $4$ can be also applied for larger multiplicities.
The case $m=5$ has become much more tiresome, and the computations are available to the interested readers upon request.
 The quasy-polynomial for the number of numerical semigroups with multiplicity $5$ and any genus has $9$ pieces, being some of them particulary large.
\section{Acknowledgement}
The first and third authors were supported by the proyects MTM2007-67433-C02-01 (Ministerio de Educaci\'on y Ciencia) and P06-FQM-01366 (Junta de Andaluc\'ia). The second author was supported by the project MTM2007-62346 (Ministerio de Educaci\'on y Ciencia) and by the research group FQM-343 (Junta de Andalucia).
The first author is also supported by the Juan de la Cierva grant JCI-2009-03896.

\end{document}